%Typed in AMSTeX

\documentstyle{amsppt}
\magnification=\magstep1
\NoBlackBoxes

 % triplenorm of #1

\topmatter

\title Isometries of $L_p$-spaces of solutions of homogeneous partial
differential equations
\endtitle

\rightheadtext{Isometries of $L_p$-spaces of solutions of differential
equations}

\author Alexander Koldobsky \endauthor
\address Division of Mathematics, Computer Science, and Statistics,
University
of Texas at San Antonio, San Antonio, TX 78249, U.S.A. \endaddress
\email koldobsk\@ringer.cs.utsa.edu \endemail

\abstract Let $ n\geq 2, A=(a_{ij})_{i,j=1}^{n}$ be a real symmetric
matrix, $a=(a_i)_{i=1}^{n}\in \Bbb R^n.$ Consider the differential
operator
$D_A = \sum_{i,j=1}^n a_{ij}{\partial^2 \over \partial x_i  \partial
x_j}+ \sum_{i=1}^n a_i{\partial \over \partial x_i}.$ Let $E$ be a
bounded domain in $\Bbb R^n,$ $p>0.$
Denote by $L_{D_A}^p(E)$ the space of solutions of the equation
$D_A f=0$ in the domain $E$ provided with the $L_p$-norm.

We prove that, for matrices $A,B,$ vectors $a,b,$ bounded domains
$E,F,$ and every $p>0$ which is not an even integer,  the space
$L_{D_A}^p(E)$
is isometric to a subspace of $L_{D_B}^p(F)$
if and only if the matrices  $A$ and $B$ have equal signatures, and the
domains
$E$ and $F$ coincide up to a natural mapping which in the most cases is
affine.
We use the extension method for $L_p$-isometries which reduces the
problem
to the question of which weighted composition operators carry solutions
of the equation $D_A f=0$ in $E$ to solutions of the equation $D_B f=0$
in $F.$
\endabstract


\endtopmatter \document \baselineskip=14pt

\head 1. Introduction \endhead

Let $n\geq 2, A=(a_{ij})_{i,j=1}^{n}$ be a real symmetric matrix,
$a=(a_i)_{i=1}^{n}\in \Bbb R^n.$ Consider the differential operator
$$D_A = \sum_{i,j=1}^n a_{ij}{\partial^2 \over \partial x_i  \partial
x_j}+ \sum_{i=1}^n a_i{\partial \over \partial x_i}$$

Let $E$ be a bounded domain in $\Bbb R^n,$ $p>0.$ Denote by
$L_{D_A}^p(E)$ the space of real functions $f\in C^2(E)$ for which $D_A
f=0$ and
$$\|f\|=(\int_E |f(x)|^p dm(x))^{1/p} < \infty$$
where $m$ is Lebesgue measure in $\Bbb R^n.$

Suppose that, for different matrices $A,B,$ vectors $a,b,$ and bounded
domains $E,F,$ the space $L_{D_A}^p(E)$ is isometric to a subspace of
$L_{D_B}^p(F).$ Does the similarity of geometric structures of the
spaces imply any equivalence of differential operators and domains ?

We shall answer this question in positive and show that, for every
$p>0$ which is not an even integer, there is a close connection between
the geometry of the space $L_{D_A}^p(E)$ and properties of $D_A$ and
$E.$

The case of the spaces of harmonic functions was considered by
A.Plotkin \cite {8}. He proved that, for $n\geq 3,$ $p\neq 2k, k\in N$
and $p\neq 2n/(n-2),$ the space $L_{\Delta}^p(E)$ is isometric to a
subspace of $L_{\Delta}^p(F)$($\Delta$ is the Laplace operator) if and
only if the domains $E$ and $F$ are similar (coincide up to the
composition of a translation, rotation, reflection and homothety). If
$p = 2n/(n-2)$ one can add an inversion to the composition. For $n=2,$
$E$ and $F$ must be similar. In \cite{5} , Plotkin's result was
extended to the case of elliptic operators $D_A$ and $D_B.$

The result of this paper generalizes Plotkin's theorem to the case of
arbitrary differential operators  $D_A$ and $D_B.$ Our main tool is the
extension method for $L_p$-isometries based on the following

\proclaim{Extension Theorem} (A.Plotkin \cite{9,10}, C.Hardin \cite{3})
Let $p>0,$ where $p$ is not an even integer,  $(X_1,\sigma_1)$ and
$(X_2,\sigma_2)$
be finite measure spaces, $Y$  a subspace of $L_p(X_1)$ containing the
constant function $1(x)\equiv 1,$ and let
T be an arbitrary linear isometry from $Y$ to $L_p(X_2).$
Then there exists a linear isometry $\tilde T: L^p (X_1, \Omega_0,
\sigma_1)
\mapsto L^p (X_2)$ such that $\tilde T \vert_Y = T.$ Here $\Omega_0$ is
the minimal $\sigma$-algebra making the functions from $Y$ measurable.
\endproclaim

By the Extension Theorem, every isometry $T: L_{D_A}^p(E_1) \mapsto
L_{D_B}^p(E_2)$ can be extended to the whole space  $L_p(E_1).$
By the classical characterizations of the isometries of $L_p$-spaces
due to S.Banach \cite{1} and J.Lamperti \cite{7}, the extension is a
weighted composition operator. Therefore, our problem can be reduced to
the following question: Which weighted composition operators carry
functions from $ L_{D_A}^p(E_1)$ to functions from $L_{D_B}^p(E_2) ?$

In \cite{6}, one can find references to other applications and
generalizations
of the extension method.

\head 2. The main result \endhead

We start with necessary definitions and notation.

Let $A=(a_{ij})_{i,j=1}^{n}, det A\neq 0$ be a real symmetric matrix.
There exists a matrix $M$ diagonalizing the matrix $A$ in the sense
that
$M^{*} A M = I_{\ell}$ for some integer $\ell, 0\leq \ell\leq n$ where
$I_{\ell} = (\ell_{ij})_{i,j=1}^n$ stands for the matrix with
$\ell_{ij}=0, i\neq j,$
$\ell_{ii}=1, 1\leq i\leq \ell$ and $\ell_{ii}=-1, \ell+1\leq i\leq n.$
We call
the number $2\ell-n$ a signature of the matrix $A.$

We denote by $D_{\ell}$ the differential operator corresponding to the
matrix
$I_{\ell}:$
$$D_{\ell} = \sum_{i=1}^{\ell} {\partial^2\over \partial x_i^2} -
\sum_{i=\ell +1}^n {\partial^2\over \partial x_i^2}.$$

For two bounded domains $E_1$ and $E_2$ in $\Bbb R^n$ with $int(cl
E_1)=E_1$
and a mapping $\tau: E_2\mapsto cl(E_1)$ of the class $C^1,$ we say
that
$E_1$ and $E_2$ coincide up to $\tau$ if
$m(E_1\setminus \tau(E_2))= 0.$

If $x\in E_2,$  $J(x)$ stands for the Jacobi matrix of $\tau$ at the
point
$x,$ and $\tau{'}(x)= det J(x)$ is the Jacobian of $\tau$ at $x.$

We say that $\tau$ is ${\ell}$-conformal at a point $x$
if $J^{*}(x) I_{\ell} J(x) = C(x) I_{\ell}$ where $C(x)\in R.$ A
mapping is
${\ell}$-conformal in a domain if it is ${\ell}$-conformal at every
point of
the domain.

For $z_0\in \Bbb R^n,$ the ${\ell}$-inversion with center $z_0$ is the
mapping
$$z\mapsto {{z-z_0}\over \|z-z_0\|_l^2} + z_0, z\in \Bbb R^n$$
where $\|z\|_{\ell}^2 =  \sum_{i=1}^{\ell} z_i^2 -
\sum_{i=\ell +1}^n z_i^2.$ A homothety with center $z_0$ and
coefficient
$t\in \Bbb R$
is the mapping $z\mapsto t(z-z_0) + z_0, z\in \Bbb R^n.$ We call
${\ell}$-similarity
a mapping which is the composition of a homothety and a mapping
preserving the metric $\|z\|_{\ell}^2$ (all such mappings are affine,
\cite{2} ).

The following characterization of ${\ell}$-conformal mappings was given
by
Liouville in 1850 for the mappings of the class $C^3.$ In 1958 the
result was extended
by Hartman \cite{4} to the $C^2$-mappings and, finally, Reshetnyak
\cite{11} formulated and proved
it without any smoothness assumptions.

\proclaim{Liouville's Theorem} Let $D$ be a domain in $\Bbb R^n, n\geq
3$ and $\tau$ a ${\ell}$-conformal mapping from $D$ to $\Bbb R^n$ where
${\ell}\in N, 1\leq {\ell}\leq n.$
Then $\tau$ is the composition of a ${\ell}$-similarity and a
${\ell}$-inversion.
\endproclaim

We are ready to formulate the main result of the paper.

\proclaim{Theorem 1} Let $n\geq 3,$ $p$ be a positive number which is
not an even integer, $E_1, E_2$ bounded domains in $\Bbb R^n$
with $int(cl E_1)=E_1,$ $A,B$ are real symmetric matrices with non-zero
determinants, $a,b\in \Bbb R^n,$ $M,N$ be the matrices diagonalizing
$A$
and $B$ and $2{\ell}-n, 2m-n$
be the signatures of $A$ and $B.$

Let $a,b\in \Bbb R^n,$ and $D_A, D_B$ be the differential operators
corresponding
to $(A,a)$ and $(B,b),$ respectively. Then:

(i) If either ${\ell}\neq m$ or ${\ell}=m$ and one of the vectors $a,b$
is zero and another
is non-zero, then  the space $ L_{D_A}^p(E_1)$ is not isometric to
a subspace of $L_{D_B}^p(E_2).$

(ii) If ${\ell}=m$ and $a=b=0$ then, for $p\neq 2n/(n-2),$ the space
$ L_{D_A}^p(E_1)$ is isometric to a subspace of $L_{D_B}^p(E_2)$
if and only if the domains $ME_1$ and $NE_2$ coincide up to a
${\ell}$-similarity $\tau.$ For $p = 2n/(n-2),$  the domains $ME_1$ and
$NE_2$ may coincide up to the composition of a ${\ell}$-similarity and
a ${\ell}$-inversion.

(iii) If ${\ell}=m$ and $a\neq 0, b\neq 0$ then the space
$ L_{D_A}^p(E_1)$ is isometric to a subspace of $L_{D_B}^p(E_2)$
if and only if the domains $ME_1$ and $NE_2$ coincide up to a
${\ell}$-similarity $\tau$ such that $JNb = |\tau{'}|^{2/n} Ma$
(since $\tau$ is an affine mapping the Jacobi matrix $J$ does not
depend on the
choice of a point.)

Finally, in all the cases where an isometric embedding $T$ exists it
has the form $Tf = \pm |det (M^{-1}JN)|^{1/p} f(M^{-1}\tau N),$
$f\in L_{D_A}^p(E_1).$
\endproclaim

In the case $n=2$ the result is different. The reason is that the
Liouville
Theorem is not valid in this case and the class of ${\ell}$-conformal
mappings is larger.

\proclaim{Theorem 2} Let $n=2$ and $p, A, B, a, b, E_1, E_2, {\ell}, m$
be as
in Theorem 1. Then:

(i) If either ${\ell}\neq m,$ or ${\ell}=m=2$ and one of the vectors
$a,b$ is zero and another
is non-zero, or ${\ell}=m=1$ and one of the numbers $\|Ma\|_1^2,
\|Nb\|_1^2$ is zero and another is non-zero, then  the space $
L_{D_A}^p(E_1)$ is not isometric to a subspace of $L_{D_B}^p(E_2).$

(ii) If ${\ell}=m$ and $a=b=0$ then, for every $p,$ the answer is the
same as  for  $p\neq 2n/(n-2)$ in Theorem 1. If ${\ell}=m=2$ and $a\neq
0, b\neq 0$ or ${\ell}=m=1$ and $\|Ma\|_1^2\neq 0, \|Nb\|_1^2\neq 0$
the answer is the same as in the part (iii) of  Theorem 1.

(iii) In the case ${\ell}=m=1$ and $\|Ma\|_1^2 = \|Nb\|_1^2 = 0$ the
class of mappings $\tau$ generating isometric embeddings is different.
The answer depends on the coordinates of the vectors
$Ma = (c_1,c_2)$ and $Nb = (d_1,d_2).$

For $c_1=\pm c_2=c\neq 0$
and $d_1=\pm d_2=d\neq 0,$ the coordinates $u_1,u_2$ of the mapping
$\tau$
are as follows:
$$\cases u_1(x_1,x_2) = (-1/pc)\ln|\gamma \exp(-pdx_1/2 \pm pdx_2/2) -
1| +
kx_1 \pm kx_2 + \alpha \\
u_2(x_1,x_2) = (\pm 1/pc)\ln|\gamma \exp(-pdx_1/2 \pm pdx_2/2) - 1| \pm
kx_1 + kx_2 + \beta \endcases$$
or
$$\cases u_1(x_1,x_2) = (-1/pd)\gamma \exp(-pdx_1/2 \pm pdx_2/2) +
(1/pc)
\ln|pcx_1 \pm pcx_2 + \delta| + \alpha \\
u_2(x_1,x_2) = (\mp 1/pd)\gamma \exp(-pdx_1/2 \pm pdx_2/2) \mp (1/pc)
\ln|pcx_1 \pm pcx_2 + \delta| + \beta \endcases$$

If  $c_1=\pm c_2=c\neq 0, d_1=\mp d_2=d\neq 0$ then
$$\cases u_1(x_1,x_2) = (-1/pd)\gamma \exp(-pdx_1/2 \mp pdx_2/2) -
(1/pc)
\ln|\mp pcx_1 + pcx_2 + \delta| + \alpha \\
u_2(x_1,x_2) = (\mp 1/pd)\gamma \exp(-pdx_1/2 \mp pdx_2/2) \pm (1/pc)
\ln|\mp pcx_1 + pcx_2 + \delta| + \beta \endcases$$
or
$$\cases u_1(x_1,x_2) = (1/pc)\ln|\gamma \exp(-pdx_1/2 \mp pdx_2/2) -
1| +
kx_1 \mp kx_2 + \alpha \\
u_2(x_1,x_2) = (\mp 1/pc)\ln|\gamma \exp(-pdx_1/2 \mp pdx_2/2) - 1| \pm
kx_1 - kx_2 + \beta \endcases$$

For  $c_1= c_2 = 0, d_1=\pm d_2=d\neq 0,$
$$\cases u_1(x_1,x_2) = (-1/pd)\gamma \exp(-pdx_1/2 \pm pdx_2/2) +
kx_1 \pm kx_2 + \alpha \\
u_2(x_1,x_2) = (\pm 1/pd)\gamma \exp(-pdx_1/2 \pm pdx_2/2) \pm
kx_1 + kx_2 + \beta \endcases$$
or
$$\cases u_1(x_1,x_2) = (-1/pd)\gamma \exp(-pdx_1/2 \pm pdx_2/2) +
kx_1 \pm kx_2 + \alpha \\
u_2(x_1,x_2) = (\mp 1/pd)\gamma \exp(-pdx_1/2 \pm pdx_2/2) \mp
kx_1 - kx_2 + \beta \endcases$$

Finally, for  $c_1=\pm c_2 = c \neq 0, d_1= d_2= 0,$
$$\cases u_1(x_1,x_2) = - (1/pc)
\ln|\mp pcx_1 + pcx_2 + \delta| + kx_1 \pm kx_2 + \alpha \\
u_2(x_1,x_2) = \pm (1/pc)
\ln|\mp pcx_1 + pcx_2 + \delta| \pm kx_1 +
kx_2 + \beta \endcases$$
or
$$\cases u_1(x_1,x_2) = \mp (1/pc)
\ln| pcx_1 \pm pcx_2 + \delta| + kx_1 \mp kx_2 + \alpha \\
u_2(x_1,x_2) = (1/pc)
\ln| pcx_1 \pm pcx_2 + \delta| \pm kx_1 - kx_2 + \beta \endcases$$

In these formulas $\alpha, \beta, \gamma, \delta, k$ are
real numbers (if $\gamma$ and $k$ are both present in a formula, one of
them must be non-zero). We use $\pm$ and $\mp$ as follows: first read
the text with
the upper signs everywhere, and then read it for the second time with
the lower signs. \endproclaim

\head 3. Weighted composition operators \endhead

We start with two facts whose simple proofs we leave to the reader.
The first one shows that, for all subspaces
of $L_p$ we are going to deal with, the $\sigma$-algebra
$\Omega_0$ appearing in the Extension Theorem is, in fact, the
$\sigma$-algebra
of all Borel sets.

\proclaim{Lemma 1} Let $E$ be a bounded open set in $\Bbb R^n,$ $H$ a
family
of continuous functions on $E$ containing the function $1(x)\equiv 1$
and
separating the points of $E$ (for every $x,y\in E,$ there exists
$f\in H$ such that $f(x)\neq f(y).)$ Then the minimal $\sigma$-algebra
of subsets of $E$ making functions from $H$ measurable is the
$\sigma$-algebra
of all Borel subsets of $E.$
\endproclaim

The second fact reduces the main question of this paper to the case
where the
matrices $A$ and $B$ are equal to $I_{\ell}$ and $I_m,$ respectively.

\proclaim{Lemma 2} Let $n\geq 2,$ $E_1, p, A, a, M, {\ell}$ be as
in Theorem 1, and
$$H = D_{\ell} + \sum_{i=1}^n \alpha_i {\partial\over\partial x_i}
\tag{1}$$
where
$\alpha=(\alpha_1,...,\alpha_n) = Ma\in R^n.$ Then the operator $T$
defined by
$$Tf(x) = |\det M|^{-1/p} f(M^{-1}x), f\in  L_{D_A}^p(E_1), x\in ME_1$$
is a linear isometry from $ L_{D_A}^p(E_1)$ onto $ L_{H}^p(ME_1).$
\endproclaim

Now we can apply the Extension Theorem to the space  $Y =
L_{H}^p(E_1).$

\proclaim{Theorem 3} Let $n\geq 2,$ $E_1, E_2$ be bounded domains in
$\Bbb R^n,$
$int(cl E_1)=E_1,$ $p>0$ and $p$ is not an even integer, ${\ell}\in N,
1\leq {\ell}\leq n,$ $\alpha\in R^n,$ and define $H$ by (1).

Then, for every isometry $T: L_{H}^p(E_1)\mapsto L^p(E_2)\cap
C^2(E_2),$
there exists a mapping $\tau: E_2\mapsto cl(E_1)$ such that:

(i) $\tau$ is of the class $C^2$ on $E_2\setminus \{x\in E_2:
T1(x)=0\}$
and $|\tau{'}(x)| = |T1(x)|^p$ on $E_2,$

(ii) $E_1$ and $E_2$ coincide up to $\tau,$

(iii) for every $f\in  L_{H}^p(E_1),$ $Tf(x) = T1(x) f(\tau(x))$ on
$E_2.$
\endproclaim

\demo{Proof} Clearly, the function $1(x)\equiv 1$ belongs to the space
$L_{H}^p(E_1).$ Besides, the space $L_{H}^p(E_1)$ separates the points
of $E_1.$ To see that, take two different points $y,z\in \Bbb R^n.$
There exists
$k$ such that $y_k\neq z_k.$ If $\alpha_k=0$ then the function
$f_k(x)\equiv
x_k$ belongs to $L_{H}^p(E_1)$ and separates $y$ and $z.$ If
$\alpha_k\neq 0$
then one of the functions $u(x)=\exp(-\alpha_k x_k)$ and
$v(x)=\exp(\alpha_k x_k)$ belongs to $L_{H}^p(E_1)$ and separates the
points.

By the Extension Theorem and Lemma 2, the isometry $T$ can be extended
to an isometry $\tilde T$ from the whole space $L^p(E_1)$ to
$L^p(E_2).$ By the classical result of J.Lamperti \cite{7},
there exists an isometric homomorphism $\phi$ from the algebra
$L^{\infty}(E_1)$ to $L^{\infty}(E_2)$ such that, for every
$f\in L^{\infty}(E_1),$ $\tilde T f = F\phi(f)$ where $F=T1.$

For the functions $f_k(x)=x_k,$ put
$\phi f_k = u_k$ and consider a mapping $\tau: E_2\mapsto \Bbb R^n$
defined by $\tau(x)=(u_1(x),...,u_n(x)).$ We are going to prove that
$\tau$ satisfies the conditions (i)-(iii).

Since $\phi$ is a homomorphism of algebras, for every polynomial
$P(x_1,...,x_n),$  $\tilde T P = FP(\tau).$ Polynomials form
a dense subset in $L_p(E_1)$ and, therefore, we have (iii).

The function $F=T1$ belongs to the class $C^2(E_2).$ Since one of the
functions $f_k,$ $u(x)=\exp(-\alpha_k x_k),$ $v(x)=\exp(\alpha_k x_k)$
belongs to
$L_{H}^p(E_1),$ one of the functions $Fu_k,$ $F\exp(-\alpha_k u_k),$
$F\exp(\alpha_k u_k)$ belongs to $C^2(E_2),$ and it follows that
$\tau$ is a mapping of the class $C^2$ on
$E_2\setminus \{x\in E_2: F(x)=0\}.$

Let us prove that $\tau(E_2)\subset cl(E_1).$ Suppose that there exists
$x_0\in E_2$ for which $\tau (x_0)\notin cl(E_1).$ Consider a
polynomial
$P(x) = A - \sum_{i=1}^n (x_i - u_i(x_0))^2$ where we choose $A>0$
so that $P$ is positive on $E_1.$ Then,
$$A = \sup_{E_2} P(\tau) = \|\phi(P)\|_{L^{\infty}(E_2)} =
\|P\|_{L^{\infty}(E_1)} = \sup_{E_1} P < A,$$
and we get a contradiction.

Let $\chi$ be the characteristic function of the set $E_1\setminus
\tau(E_2).$
Then $\tilde T \chi = F\chi(\tau) = 0.$ Since $\tilde T$ is an isometry
we get $\chi=0$ which means that $m(E_1\setminus \tau(E_2)) = 0,$
and the domains $E_1, E_2$ coincide up to $\tau.$ We have proved (ii).

To finish the proof of (i), note that, for every function $f\in
L_p(E_1),$
$$\|f\|_{L_p(E_1)}^p = \int_{E_1} |f(y)|^p dm(y) = $$
$$\int_{E_2} |f(\tau(x))|^p |\tau^{'}(x)| dm(x) =$$
$$\|\tilde T f\|_{L_p(E_2)}^p = \int_{E_2} |f(\tau(x))|^p |F(x)|^p
dm(x).$$
(We made a change of variables $y=\tau(x).$) Since $f$ is an arbitrary
function and  $m(E_1\setminus \tau(E_2)) = 0$ we get $|\tau{'}|=|F|^p,$
which completes the proof of the theorem.
\qed \enddemo

Part (iii) of Theorem 3 shows that every isometry from $L_{H}^p(E_1)$
to $L^p(E_2)\cap C^2(E_2)$ is generated by a mapping $\tau.$ Now we are
going to choose those mappings for which the images of functions
from $L_{H}^p(E_1)$ are solutions of another differential equation.

We need the following elementary fact.

\proclaim{Lemma 3} Let $n,{\ell}\in N, 1\leq {\ell}\leq n.$ Put
$\epsilon_i = 1,
1\leq i \leq {\ell}$ and $\epsilon_i = -1, {\ell}< i \leq n.$ Suppose
that a real symmetric matrix
$B = (b_{ij})_{i,j=1}^n,$ a vector $a\in R^n,$ and a number $c\in R$
satisfy the following: for any choice of complex numbers $s_1,...,s_n,$
the equality $\sum_{i=1}^n \epsilon_i s_i^2 + a_is_i =0$ implies
$\sum_{i,j=1}^n  b_{ij}s_is_j +  \sum_{i=1}^n ca_is_i = 0.$

Then the matrices $B$ and $I_{\ell}$ differ by a constant multiple
only.
\endproclaim

\proclaim{Theorem 4}  Let $n\geq 2,$ $E_1, E_2$ be bounded domains in
$\Bbb R^n,$
$\tau: E_2\mapsto cl(E_1)$ a mapping of the class $C^2,$ $\alpha\in
\Bbb R^n,$
and $H$ the differential operator defined by (1).

Consider any real functions $b_{ij}(x), b_i(x), i,j=1,...,n$ on $E_2$
and denote by $D_B$ the differential operator
$$D_B =  \sum_{i,j=1}^n b_{ij}(x){\partial^2 \over \partial x_i
\partial x_j}+ \sum_{i=1}^n b_i(x){\partial \over \partial x_i}.$$

Suppose that there exists a function $F\in L_{D_B}(E_2)$ such that, for
every
$f\in L_H^{\infty} (E_1),$ the function $Ff(\tau)$ belongs to
$L_{D_B}(E_2).$
Then:

(i) for every $x\in E_2\setminus \{x\in E_2: F(x)=0\},$ the matrix
$B(x)=(b_{ij}(x))_{i,j=1}^n$ has the signature $2l-n.$

(ii) there exists a real function $C: E_2\mapsto \Bbb R$ such that, for
every
$x\in E_2\setminus \{x\in E_2: F(x)=0\},$
$$J^{*}(x) B(x) J(x) = {C(x)\over F(x)} I_{\ell}.$$
\endproclaim

\demo{Proof} Let $\tau = (u_1,...,u_n).$ For any $a\in R^n,$
 $$\sum_{k=1}^n a_k x_k \in L_H^{\infty} (E_1) \Longleftrightarrow
\sum_{k=1}^n \alpha_k a_k = 0. \tag{2}$$
On the other hand, if $\sum_{k=1}^n a_k x_k \in L_H^{\infty} (E_1)$
then $F\sum_{k=1}^n a_k u_k \in L_{D_B}(E_2).$
Since $F\in L_{D_B}(E_2)$ the latter condition gives
$$\sum_{k=1}^n a_k \bigl(\sum_{i,j=1}^n b_{ij}(x)(2{\partial F\over
\partial x_i}{\partial u_k\over \partial x_j} + F{\partial^2 u_k\over
\partial x_i
\partial x_j}) +
\sum_{i=1}^n b_i(x) F{\partial u_k \over \partial x_i}\bigr) = 0.
\tag{3}$$

Since (2) implies (3) for every vector $a,$ the
coefficients at $a_k$'s must be proportional. It means that , for every
$x\in E_2,$ there exists $c(x)\in \Bbb R$ such that, for each
$k=1,...,n,$

$$\sum_{i,j=1}^n b_{ij}(x)(2{\partial F\over \partial x_i}(x){\partial
u_k\over \partial x_j}(x) + F(x){\partial^2 u_k\over \partial x_i
\partial x_j}(x)) +
\sum_{i=1}^n b_i(x) F(x){\partial u_k \over \partial x_i}(x) = c(x)
\alpha_k. \tag{4}$$

Consider the function $\exp(x,s)$ where $s=(s_1,...s_n)$ is a $n$-tuple
of complex numbers and $(x,s)$ stands for the scalar product. Clearly,
$H(\exp(x,s))=0$ if and only if
$$\sum_{i=1}^n \epsilon_i s_i^2 + \alpha_is_i =0 \tag{5}$$
where the numbers $\epsilon_i$ are the same as in Lemma 3.

On the other hand, $D_B(F\exp(\tau(x),s)) = 0$ and using (4) and the
fact that $F\in L_{D_B}(E_2)$  we get

$$F(x)\sum_{k,m=1}^n s_k s_m \sum_{i,j=1}^n b_{ij}(x)
{\partial u_k \over \partial x_i}(x){\partial u_m \over \partial
x_j}(x) +
c(x)\sum_{k=1}^n \alpha_k s_k = 0. \tag{6}$$

Thus, for any choice of complex numbers $s_1,...,s_n,$ (5) implies
(6).
It means that, for every $x\in  E_2\setminus \{x\in E_2: F(x)=0\},$
the matrix $J^{*}(x) B(x) J(x)$ satisfies the conditions of Lemma 3.
Therefore, there exists a function $C:E_2\mapsto R$ such that
$$J^{*}(x) B(x) J(x) = {C(x)\over F(x)} I_{\ell},$$
and we get (ii). Part (i) follows from the uniqueness of the
diagonalization.
\qed \enddemo

\head 4. Proof of the main result \endhead

We are ready to prove Theorems 1 and 2. The first part of the proof
applies to both of the cases $n>2$ and $n=2.$

Using Lemma 2 one can reduce the problem to the case of diagonal
matrices.
Let $T$ be an isometry from  $ L_{D_A}^p(E_1)$ to $L_{D_B}^p(E_2)$ and
define differential operators $H$ and $G$ by
$$H = D_{\ell} + \sum_{i=1}^n c_i {\partial\over\partial x_i},
 G = D_{m} + \sum_{i=1}^n d_i {\partial\over\partial x_i}$$
where $c=Ma, d=Nb.$

By Lemma 2, the operators

$$T_1f(x)= |\det M|^{-1/p} f(M^{-1}x),
T_2g(x)= |\det N|^{-1/p} g(N^{-1}x)$$
are isometries from $ L_{D_A}^p(E_1)$ and $L_{D_B}^p(E_2)$ onto
$ L_{H}^p(ME_1)$ and
$L_{G}^p(NE_2),$ respectively. Therefore, $S=T_2 T T_1^{-1}$ is an
isometry from $L_{H}^p(ME_1)$ to $L_{G}^p(NE_2).$

Put $F=S1$ and $E=\{x\in NE_2: F(x)=0\}.$ By Theorem 3, there exists
a mapping $\tau: NE_2\mapsto ME_1$ of the class $C^2$ on
$NE_2\setminus E$ such that the domains $ME_1$ and $NE_2$ coincide
up to $\tau,$ $|\tau'|\equiv |F|^p$ and, for every $f\in
L_{H}^p(ME_1),$
$Sf = Ff(\tau).$

Clearly, the mapping $\tau$ satisfies the conditions of Theorem 4 with
$B(x)=G$ for every $x,$
so the matrices $H$ and $G$ have equal signatures which means
that $\ell = m.$ Besides, there exists a real function
$C: NE_2\mapsto \Bbb R$ such that
$$J^{*}(x) I_{\ell} J(x) = {C(x)\over F(x)} I_{\ell} \tag{7}$$
for every $x\in NE_2\setminus E.$

Thus, $\tau$ is a $\ell$-conformal mapping on $NE_2\setminus E.$

Calculating the determinants in both sides of (7) we get
$$J^{*}(x) I_{\ell} J(x) = |\tau'(x)|^{2/n} I_{\ell} \tag{8}$$

Let $\epsilon_i = 1, 1\leq i \leq {\ell}$ and $\epsilon_i = -1,
{\ell}< i \leq n.$ If $\tau = (u_1,...,u_n)$ then, for every function
$f\in L_H^p(ME_1),$
$$\split 0=G(Ff(\tau)) = \sum_{k=1}^n {\partial f\over \partial
x_k}(\tau)
\sum_{i=1}^n \biggl((\epsilon_i (2{\partial F\over \partial
x_i}{\partial u_k\over \partial x_i} +
F{\partial^2u_k\over \partial x_i^2}) + \\
F{\partial u_k\over \partial x_i} d_i - c_k F|\tau'(x)|^{2/n}\biggr)
\endsplit \tag{9}$$
We used (8) and the fact that $F\in L_G^p(NE_2).$

Starting from this point we consider the cases $n\geq 3$ and $n=2$
separately.

\subhead{The case \bf$n\geq 3$} \endsubhead
It follows from (8) that $\tau$ is a $\ell$-conformal mapping on
$NE_2\setminus E.$ By Liouville's theorem, the mapping $\tau$
is either a $\ell$-similarity or the composition of a $\ell$-similarity
and a $\ell$-inversion
on every connected subset $U$ of $NE_2\setminus E.$ The Jacobian
of the $\ell$-inversion with center $x_0$ is equal to
$\|x-x_0\|_{\ell}^{-2n}$ (see \cite{2}) and the Jacobian of any
$\ell$-similarity is a constant, so we have
$\tau'(x) = k\|x-x_0\|_{\ell}^{-2n}$ on $U.$

By Theorem 4,  $|\tau'|\equiv |F|^p,$ and
we have $|F| = k^{1/p}\|x-x_0\|_{\ell}^{-2n/p}$ on $U.$
Clearly, $k\neq 0$ because $U\cap E =\emptyset.$ Therefore,
there exists a constant $\alpha > 0$ such that
$|F(x)|>\alpha$ on $U.$ The function $F$ is continuous, so
$cl (U) \cap E = \emptyset$ for every connected subset $U$ of
$NE_2\setminus E.$ This is possible only if $E=\emptyset.$
Thus, the mapping $\tau$ is either a $\ell$-similarity or the
composition of a $\ell$-similarity and a $\ell$-inversion
on the whole set $NE_2.$

Suppose that $\tau$ is the composition of a $\ell$-similarity and a
$\ell$-inversion on $NE_2.$ Since the function $F(x) =
k^{1/p}\|x-x_0\|_{\ell}^{-2n/p}$ belongs to the space $L_{G}^p(NE_2)$
we have
$$\split 0=G(F)=(-nk/p) \sum_{i=1}^n d_i \|x-x_0\|_{\ell}^{-2n/p - 2}
2(x_i-
(x_0)_i) \epsilon_i  + \\
(nk/p)(n/p+1) \sum_{i=1}^n  \|x-x_0\|_{\ell}^{-2n/p - 4}
4(x_i-(x_0)_i)^2
\epsilon_i  + \\
(2kn^2/p) \|x-x_0\|_{\ell}^{-2n/p - 2} \endsplit \tag{10}$$
for every $x\in NE_2.$ It is easy to see that (10) implies $p=2n/(n-2)$
and $d=0.$

Simple calculations show that, for $F(x)= k\|x-x_0\|_{\ell}^{2-n}$ and
$u_m(x)=(x_0)_m + (x_m - (x_0)_m)/\|x-x_0\|_{\ell}^{2}$ (these are the
coordinate functions for $\ell$-inversion),
$$  \sum_{k=1}^n \epsilon_i (2{\partial F\over \partial x_i}
{\partial u_k\over \partial x_i} +
F{\partial^2u_k\over \partial x_i^2}) \equiv 0$$
Now (9) implies $c=0.$

Thus, the mapping $\tau$ can contain a $\ell$-inversion only if
$p=2n/(n-2)$ and $c=d=0.$ On the other hand, if $p=2n/(n-2),$ $c=d=0$
and $\tau$ is the composition of a $\ell$-similarity and a
$\ell$-inversion
then, by (10), $|\tau'|^{1/p}\in L_{G}^p(NE_2)$ and , by (9),
$|\tau'|^{1/p}f(\tau)\in  L_{G}^p(NE_2)$ for every $f\in L_H^p(ME_1).$
Therefore, $Sf=|\tau'|^{1/p} f(\tau)$ is an isometry from $L_H^p(ME_1)$
to $L_{G}^p(NE_2).$ We have proved part (ii) of the theorem.

If $p\neq 2n/(n-2)$ or one of the vectors $c, d$ is non-zero the
mapping
$\tau$ is a $\ell$-similarity. Therefore, $F = |\tau'|^{1/p}$ is a
constant
function, and the coordinate functions $u_k$ of the mapping $\tau$
are affine. The equality (9) implies
$$\sum_{i=1}^n F{\partial u_k\over \partial x_i}d_i =
F|\tau'|^{2/n}c_k$$
for every $k=1,...,n.$ It means that $Jd=|\tau'|^{2/n}c$ which proves
(iii)
and, besides, shows that an isometric embedding does not exist if one
of the vectors $c,d$ is zero and another is non-zero. We have proved
Theorem 1.

\subhead{The case \bf$n=2$} \endsubhead
In this case, Liouville's theorem is no longer valid. We consider
the cases $\ell=2$ and $\ell=1$ separately.

First, let $\ell=2.$ It follows from (8) that the mapping $\tau$
is either holomorphic or antiholomorphic, and
$${\partial u_1\over \partial x_1}=\pm {\partial u_2\over \partial
x_2},
{\partial u_1\over \partial x_2}=\mp{\partial u_2\over \partial x_1}
\tag{11}$$
(Read this formula and the following text with the upper signs first,
and then
read it for the second time with the lower signs.)

The equality (9) shows that, for $k=1,2,$
$$(2{\partial F\over \partial x_1} + d_1F){\partial u_k\over \partial
x_1} +
(2{\partial F\over \partial x_2} + d_2F){\partial u_k\over \partial
x_2} =
F|\tau'|c_k \tag{12}$$
We get from (11) and (12) that , for $k=1,2,$
$${{\partial F\over \partial x_k}\over F} = \pm (1/2)
(c_1{\partial u_1\over \partial x_k} + c_2{\partial u_2\over \partial
x_k}-d_k) \tag{13}$$

Solving this system of equations with respect to $F$ we get
$$F(x) = \exp(\pm (1/2)(c_1u_1(x)+c_2u_2(x)-d_1x_1-d_2x_2)) \tag{14}$$

Since $F\in L_G^p(NE_2)$ we get using (11), (13) and (14) that
$$\split G(F) = (1/4)F((c_1{\partial u_1\over \partial x_1} +
c_2{\partial u_2\over \partial x_1}-d_1)^2 +
(c_1{\partial u_1\over \partial x_2} + c_2{\partial u_2\over \partial
x_2}-d_2)^2 + \\
 2(c_1{\partial^2 u_1\over \partial x_1^2} +
c_2{\partial^2 u_2\over \partial x_1^2}) \pm
2(c_1{\partial^2 u_1\over \partial x_2^2} +
c_2{\partial^2 u_2\over \partial x_2^2})+ \\
2d_1(c_1{\partial u_1\over \partial x_1} + c_2{\partial u_2\over
\partial x_1}-d_1) + 2d_2(c_1{\partial u_1\over \partial x_2} +
c_2{\partial u_2\over \partial x_2}-d_2) = \\
(1/4)F(\pm (c_1^2+c_2^2)|\tau'| - (d_1^2+d_2^2)) = 0 \endsplit
\tag{15}$$

If $c\neq 0, d\neq 0$ we get from the latter equality that
$|\tau'|$ is constant.(Note that, by (13), $F$ is non-zero at every
point.) Clearly, $\tau$ is a similarity, and we get
$Jd=|\tau'|c$ in the same way as in the case $n\geq 3.$
If one of the vectors $c,d$ is zero and another is non-zero,
(15) is impossible, so an isometry does not exist.
If $c=d=0,$ (9) implies ${\partial F\over \partial x_1}\equiv 0$ and
${\partial F\over \partial x_2}\equiv 0.$ Therefore, $|\tau'|=|F|^p$
is constant and $\tau$ is a similarity.
This finishes the proof in the case $\ell=2.$

Let $\ell = 1.$ Then, instead of (11), we get
$${\partial u_1\over \partial x_1}=\pm {\partial u_2\over \partial
x_2},
{\partial u_1\over \partial x_2}=\pm{\partial u_2\over \partial x_1}
\tag{16}$$
and $\tau' = \pm (({{\partial u_1}\over{ \partial x_1}})^2 -
({{\partial u_1}\over {\partial x_2}})^2). $

We get from (9) that, for $k=1,2,$
$$(2{\partial F\over \partial x_1} + d_1F){\partial u_k\over \partial
x_1} +
(-2{\partial F\over \partial x_2} + d_2F){\partial u_k\over \partial
x_2} =
F|\tau'|c_k \tag{17}$$
It follows from (16) and (17) that
$${{\partial F\over \partial x_k}\over F} = \pm (1/2)
(c_1{\partial u_1\over \partial x_k} - c_2
{\partial u_2\over \partial x_k}-(-1)^k d_k) \tag{18}$$
$k=1,2,$ and we can calculate $F:$
$$F(x) = \exp(\pm (1/2)(c_1u_1(x)-c_2u_2(x)-d_1x_1+d_2x_2)) \tag{19}$$

Similarly to the case $\ell = 2$ we get
$\pm (c_1^2-c_2^2)|\tau'| - (d_1^2-d_2^2)) = 0.$ We can finish the
proof
in the same way as for $\ell = 2$ if either $\|c\|_1^2\neq 0,
\|d\|_1^2\neq 0$ or one of the numbers  $\|c\|_1^2, \|d\|_1^2$
is zero and another is non-zero.

Finally, consider the case where $c_1^2=c_2^2$ and $d_1^2=d_2^2.$
Suppose we have (16) with $+$. Since $|F|=|\tau'|^{1/p}= (({{\partial
u_1}\over{ \partial x_1}})^2 - ({{\partial u_1}\over {\partial
x_2}})^2)^{1/p}$ we can
write (18) in the following form:
$$\cases {\partial u_1\over \partial x_1}{\partial^2 u_1\over \partial
x_1^2} -
{\partial u_1\over \partial x_2}{\partial^2 u_1\over {\partial
x_1\partial x_2}}
= (p/4) (({\partial u_1\over \partial x_1})^2 -
({\partial u_1\over \partial x_2})^2) (c_1{\partial u_1\over \partial
x_1} -
c_2{\partial u_1\over \partial x_2} - d_1) \\
{\partial u_1\over \partial x_1}{\partial^2 u_1\over {\partial
x_1\partial x_2}} - {\partial u_1\over \partial x_2}{\partial^2
u_1\over {\partial x_1^2}}
= (p/4) (({\partial u_1\over \partial x_1})^2 -
({\partial u_1\over \partial x_2})^2) (c_1{\partial u_1\over \partial
x_2} -
c_2{\partial u_1\over \partial x_1} + d_2) \endcases \tag{20}$$
Adding and subtracting the equations (20) we get
$$\cases {\partial^2 u_1\over \partial x_1^2} + {\partial^2 u_1\over
{\partial x_1\partial x_2}} = (p/4)({\partial u_1\over \partial x_1}+
{\partial u_1\over \partial x_2})((c_1-c_2)({\partial u_1\over \partial
x_1}+
{\partial u_1\over \partial x_2}) - (d_1-d_2)) \\
{\partial^2 u_1\over \partial x_1^2} - {\partial^2 u_1\over {\partial
x_1\partial x_2}} = (p/4)({\partial u_1\over \partial x_1}-
{\partial u_1\over \partial x_2})((c_1+c_2)({\partial u_1\over \partial
x_1}-
{\partial u_1\over \partial x_2}) - (d_1+d_2)) \endcases  \tag{21}$$

If we found a function $u_1$ satisfying (21) and defined $u_2$
so that (16) holds then the mapping $\tau=(u_1,u_2)$ would generate
an isometry from $L_H^p(ME_1)$ to  $L_G^p(NE_2)$  because
$|\tau'|^{1/p}\in
L_G^p(NE_2)$ and $|\tau'|^{1/p}f(\tau)\in L_G^p(NE_2)$ for every
$f\in L_H^p(ME_1).$

First, assume that $c_1=c_2=c\neq 0$ and $d_1=d_2=d\neq 0.$ Then (21)
implies
$$\cases {\partial^2 u_1\over \partial x_1^2} + {\partial^2 u_1\over
{\partial x_1\partial x_2}} = {\partial\over {\partial x_1}}({\partial
u_1\over \partial x_1}+{\partial u_1\over \partial x_2}) = 0 \\
{\partial\over {\partial x_1}}({\partial u_1\over \partial
x_1}-{\partial u_1\over \partial x_2}) = (p/2)({\partial u_1\over
\partial x_1}-{\partial u_1\over \partial x_2})(c({\partial u_1\over
\partial x_1}-{\partial u_1\over \partial x_2})-d) \endcases \tag{22}$$
Modify (22) using (16) to get
$$\cases {\partial^2 u_1\over \partial x_1^2} + {\partial^2 u_1\over
{\partial x_1\partial x_2}} = {\partial\over {\partial x_2}}({\partial
u_1\over \partial x_1}+{\partial u_1\over \partial x_2}) = 0 \\
-{\partial\over {\partial x_1}}({\partial u_1\over \partial
x_1}-{\partial u_1\over \partial x_2}) = (p/2)({\partial u_1\over
\partial x_1}-{\partial u_1\over \partial x_2})(c({\partial u_1\over
\partial x_1}-{\partial u_1\over \partial x_2})-d) \endcases \tag{23}$$

It follows from (22) and (23) that
$${\partial u_1\over \partial x_1}+{\partial u_1\over \partial x_2} =
K=const\tag{24}$$
Integrating the second equalities in (22) and (23) and using (16) we
get
$${\partial u_1\over \partial x_1}-{\partial u_1\over \partial x_2} =
{{d\gamma \exp((-pd/2)x_1 + (pd/2)x_2)}\over {c(\gamma \exp((-pd/2)x_1
+ (pd/2)x_2) - 1)}} \tag{25}$$
for some $\gamma\in \Bbb R.$ By (24) and (25),
$$u_1 = (-1/cp) \ln |\gamma \exp((-pd/2)x_1 + (pd/2)x_2) - 1| +
(K/2)x_1
+ (K/2)x_2 + \alpha$$
where $\alpha\in \Bbb R.$ Now we can use (16) to find $u_2:$
$$u_2 = (1/cp) \ln |\gamma \exp((-pd/2)x_1 + (pd/2)x_2) - 1| + (K/2)x_1
+ (K/2)x_2 + \beta.$$

Similarly, one can calculate $u_1$ and $u_2$ for all other cases
considered in the third part of Theorem 2. Note that, in every case,
one gets the first solution if (16) holds with positive signs and the
second solution appears if we have (16) with negative signs.
We have proved Theorem 2.

\enddocument